\input amstex
\documentstyle{amsppt}
\magnification=1200
\hcorrection{0.1in}
\vcorrection{-0.4in}
\define\volume{\operatorname{vol}}
\define\op#1{\operatorname{#1}}
\define\svolball#1#2{{\volume(\underline B_{#2}^{#1})}}
\define\svolann#1#2{{\volume(\underline A_{#2}^{#1})}}
\define\sball#1#2{{\underline B_{#2}^{#1}}}
\define\svolsp#1#2{{\volume(\partial \underline B_{#2}^{#1})}}

\NoRunningHeads
\NoBlackBoxes
\topmatter
\title Quantitative Volume Space From Rigidity with lower Ricci curvature bound II
\endtitle
\author Lina Chen, \footnote{Supported partially by a research fund from Capital Normal University.
\hfill{$\,$}} Xiaochun Rong \footnote{Supported partially by NSF
Grant DMS 0203164 and by a research fund from Capital Normal University.
\hfill{$\,$}} and Shicheng Xu\footnote{Supported partially by NSFC Grant 11401398
and by a reach fund from Capital Normal University. \hfill{$\,$}}
\endauthor
\address Mathematics Department, Capital Normal University, Beijing,
P.R.C.
\endaddress
\email chenlina\_mail\@163.com, \,\, shichxu\@foxmail.com
\endemail
\address Mathematics Department, Rutgers University
New Brunswick, NJ 08903 USA
\endaddress
\email rong\@math.rutgers.edu
\endemail

\abstract This is the second paper of two in a series under the same title 
([CRX]); both study the quantitative volume space form rigidity conjecture: 
a closed $n$-manifold of Ricci curvature at least $(n-1)H$, $H=\pm 1$ or $0$ 
is diffeomorphic to a $H$-space form if for every ball of definite size on 
$M$, the lifting ball on the Riemannian universal covering space of the ball 
achieves an almost maximal volume, provided the diameter of $M$ is bounded 
for $H\ne 1$.

In [CRX], we verified the conjecture for the case that $M$ or its Riemannian 
universal covering space $\tilde M$ is not collapsed for $H=1$ or $H\ne 1$ 
respectively. In the present paper, we will verify this conjecture for the 
case that Ricci curvature is also bounded above, while the above 
non-collapsing condition is not required.
\endabstract
\endtopmatter
\document

\head 0. Introduction
\endhead

\vskip4mm

This is the second paper of two in a series under the same title, concerning 
the quantitative version of the following volume space form rigidity.

Let $M$ be a compact $n$-manifold of Ricci curvature bounded below by 
$(n-1)H$, a constant. For $p\in M$ and $r>0$, the volume of the $r$-ball at 
$p$, $\op{vol}(B_r(p))\le \svolball{H}{r}$, and ``$=$'' if and only if the 
open ball $B_r(p)$ is isometric to $\sball{H}{r}$ (Bishop volume 
comparison), which denotes the $r$-ball in the $n$-dimensional simply 
connected $H$-space form.

The following statement is a consequence of the Bishop volume comparison.

\proclaim{Theorem 0.1} (Volume space form rigidity) Let $\rho>0$. If a 
compact $n$-manifold $M$ satisfies $$\op{Ric}_M\ge (n-1)H,\quad 
\frac{\op{vol}(B_\rho(x^*))} {\svolball{H}{\rho}}=1,\quad \forall \, x\in 
M,$$ then $M$ is isometric to a space form of constant curvature $H$,
where $\pi^*: (\widetilde{B_\rho(x)},x^*)\to (B_\rho(x),x)$ is the 
(incomplete) Riemannian universal covering space. 
\endproclaim

All $H$-space forms satisfy the local volume condition in Theorem 0.1. On 
the other hand, given any $\rho, \epsilon>0$ and $H=\pm 1$ or $0$, there is 
a $H$-space form which contains a point $x$ such that $\op{vol}(B_\rho(x)) 
<\epsilon$ i.e., $B_\rho(x)$ is collapsed.

In [CRX], we proposed the following quantitative version of Theorem 0.1.

\example{Conjecture 0.2} (Quantitative volume space form rigidity)
Given $n, \rho,d>0$ and $H=\pm 1$ or $0$, there exists a constant 
$\epsilon(n,\rho,d)>0$ such that for any $0<\epsilon<\epsilon(n,\rho,d)$, if 
a compact $n$-manifold $M$ satisfies $$\op{Ric}_M\ge (n-1)H,\quad d\ge 
\op{diam}(M),\quad \frac{\volume (B_\rho(x^*))}{\svolball{H}{\rho}}\ge 
1-\epsilon,\quad\forall\, x\in M,$$ then $M$ is diffeomorphic and 
$\Psi(\epsilon|n,\rho,d)$-close in the Gromov-Hausdorff topology
to a space form of constant curvature $H$, where $d=\pi$ or $1$ when $H=1$ 
or $0$ respectively, where $\Psi(\epsilon|n,\rho,d)\to 0$ as $\epsilon\to 0$ 
while $n, \rho$ and $d$ are fixed.
\endexample

Note that Conjecture 0.2 for $H=-1$ does not hold if one removes a bound on 
diameter (see [CRX]). On the other hand, for $H\ne -1$, $M$ in Conjecture 
0.2 may have arbitrarily small volume i.e., $M$ is collapsed.

By the volume convergence ([Co2]), Conjecture 0.2 implies the following:

\example{Conjecture 0.3} (Non-collapsing on Riemannian universal cover) 
Given $n, \rho, d>0$, $H=\pm 1$ or $0$, there exist constants, 
$\epsilon(n,\rho,d), v(n,\rho,d)>0$, such that if a compact $n$-manifold $M$ 
satisfies
$$\op{Ric}_M\ge (n-1)H,\quad d\ge \op{diam}(M),\quad\frac{\volume
(B_\rho(x^*))}{\svolball{H}{\rho}}\ge 1-\epsilon(n,\rho,d),\quad\forall\, 
x\in M,$$ then $\op{vol}(B_1(\tilde x))\ge v(n,\rho,d)>0$, where $\tilde x$ 
is a point in the Riemannian universal covering space.
\endexample

In [CRX], among other things we proved that Conjecture 0.3 implies 
Conjecture 0.2 for $H\ne 1$, and for $H=1$, Conjecture 0.2 holds when $M$ is 
not collapsed. Precisely, the following theorem is a combination of Theorem 
A, B and C in [CRX] (corresponding to $H=1, -1$ and $0$).

\proclaim{Theorem 0.4} Given $n,\rho, d, v>0$ and $H=\pm 1$ or $0$, there 
exists a constant $\epsilon(n,\rho,d,v)>0$ such that for any $0<\epsilon< 
\epsilon(n,\rho,d,v)$, if a compact $n$-manifold $M$ satisfies 
$$\op{Ric}_M\ge (n-1)H,\,d\ge \op{diam}(M),\,\op{vol}(B_1(z_0))\ge v,\, 
\frac{\op{vol}(B_\rho(x^*))}{\svolball{H}{\rho}}\ge 1-\epsilon,\, \forall\, 
x\in M,$$
then $M$ is diffeomorphic and $\Psi(\epsilon|n,\rho,d,v)$-close to a space form of constant curvature $H$,
where $d=\pi$ or $1$ when $H=1$ or $0$, $z_0\in M$ or $z_0\in \tilde M$ when 
$H=1$ or $H\ne 1$.
\endproclaim

For $H=1$, Theorem 0.4 generalizes the differential sphere theorem in [CC2] 
(cf. [Pe], [Co1], see Remark 0.7 in [CRX]), and for $H=-1$, Theorem 0.4 is 
equivalent to a quantitative version of the maximal volume entropy rigidity 
in [LW] (see Theorem D, Corollary 0.6 in [CRX]).

In the present paper, we will verify Conjecture 0.2 under an additional 
assumption: Ricci curvature is also bounded above (Theorem D). This 
regularity condition allows us to find a nearby metric of almost constant 
sectional curvature (Theorem B) by smoothing method ([DWY]) via renormalized 
Ricci flows in sense of [TW]. As an application we verify Conjecture 0.3 in 
this case (Theorem C).

We now begin to state the main results in this paper.

The first result says that under bounded Ricci curvature, the almost 
maximality of volume on local coverings measures how far the metric from 
being an $H$-Einstein metric (compare to Remark 0.5).

\proclaim{Theorem A} Given $n,\rho, \Lambda>0$, $H=\pm 1$ and $0$, there 
exists a constant, $\epsilon(n,\rho,\Lambda)>0$, such that for 
$0<\epsilon<\epsilon(n,\rho,\Lambda)$, if a compact Riemannian $n$-manifold 
$(M,g)$ satisfies $$\Lambda\ge \op{Ric}(g)\ge (n-1)H,\quad \frac{\volume
(B_\rho(x^*))}{\svolball{H}{\rho}}\ge 1-\epsilon,\quad\forall\, x\in M,$$
then $g$ is almost Einstein in $L^p$-sense for any $p\ge 1$ i.e., 
$$-\kern-1em\int_M|\op{Ric}(g)-(n-1)Hg|^p<\Psi(\epsilon|n,\rho,\Lambda,p).$$ 
\endproclaim

The additional upper bound on Ricci curvature implies a uniform 
$C^{1,\alpha}$-Harmonic radius on $B_{\frac \rho2}(x^*)$ (Lemma 1.3), and a 
local version of Theorem A on $B_{\frac \rho2}(x^*)$ (Lemma 1.4). By a 
packing argument via relative volume comparison, we obtain Theorem A.

Consider the Ricci flow on $(M,g)$; following [DWY] we see that bounded 
Ricci curvature and a uniform $C^{1,\alpha}$-Harmonic radius on $B_{\frac 
\rho2}(x^*)$ (independent of $x$) imply that the Ricci flow on $M$ exists 
for a definite time (Theorem 1.5), and that the renormalized Ricci flow 
([TW]) preserves the almost Einstein property in $L^p$-sense (Lemma 1.7). 
Using the two properties, we will prove the following strong smoothing 
result.

\proclaim{Theorem B} (Smoothing to almost constant curvature)  Given $n,
\rho, \Lambda, \delta>0$ and $H=\pm 1$ or $0$, there exists a constant, 
$\epsilon(n,\rho,\Lambda,\delta)>0$, such that for any 
$0<\epsilon<\epsilon(n,\rho,\Lambda,\delta)$, if a compact $n$-manifold 
$(M,g)$ satisfies $$\Lambda\ge \op{Ric}(g)\ge (n-1)H,\quad \frac{\volume
(B_\rho(x^*))}{\svolball{H}{\rho}}\ge 1-\epsilon,\quad\forall\, x\in M,$$
then $M$ admits a metric $g'$ such that $|g'-g|<\delta$ and for any $0\le 
k<\infty$, $$|\op{Rm}(g')|_{C^k,M}\le C(n,\rho,\Lambda,\delta,k),\,
H-\Psi(\delta,\epsilon|n,\rho,\Lambda)\le \op{sec}_{g'}\le 
H+\Psi(\delta,\epsilon|n,\rho,\Lambda).$$
\endproclaim

Using the existence of a nearby metric of almost constant sectional 
curvature, we are able to verify Conjecture 0.3 for the case bounded Ricci 
curvature.

\proclaim{Theorem C} Given $n, \rho, d, \Lambda>0$ and $H=\pm 1$ or $0$, 
there exist positive constants, $\epsilon(n,\rho,d,\Lambda), 
v(n,\rho,d,\Lambda)>0$, such that if a compact $n$-manifold $M$ satisfies
$$\Lambda\ge \op{Ric}_M\ge (n-1)H,\, d\ge \op{diam}(M),\,\frac{\volume
(B_\rho(x^*))}{\svolball{H}{\rho}}\ge 1-\epsilon(n,\rho,d,\Lambda),\quad\forall\, x\in M,$$
then $\tilde M$ is not collapsed i.e., $\op{vol}(B_1(\tilde p))\ge v(n,\rho,d,\Lambda)$
for any $\tilde p\in \tilde M$, where $d=\pi$ or $\infty$ when $H=1$ or $-1$.
\endproclaim

By Theorem C, we can apply Theorem 0.4 for $H\ne 1$ to verify Conjecture 0.2 
for the case of bounded Ricci curvature. For $H=1$ (bounded Ricci 
curvature), by the higher regularity of a nearby metric in Theorem B
we are able to strengthen Theorem 3.5 in [CRX] to conclude Conjecture 0.2 in 
this case.

\proclaim{Theorem D} Given $n, \rho, d,\Lambda>0$ and $H=\pm 1$ or $0$, 
there exists a constant $\epsilon(n,\rho,d,\Lambda)>0$ such that for any 
$0<\epsilon<\epsilon(n,\rho,d,\Lambda)$, if a compact $n$-manifold $M$ 
satisfies
$$\Lambda\ge \op{Ric}_M\ge (n-1)H,\quad d\ge \op{diam}(M),\quad \frac{\volume
(B_\rho(x^*))}{\svolball{H}{\rho}}\ge 1-\epsilon,\quad\forall\, x\in M,$$
then $M$ is diffeomorphic and $\Psi(\epsilon|n,\rho,d,\Lambda)$-close to a 
space form of constant curvature $H$, where $d=\pi$ or $1$ when $H=1$ or $0$ 
respectively.
\endproclaim

As mentioned in the above, there is a uniform lower bound on 
$C^{1,\alpha}$-harmonic radius on $B_{\frac \rho2}(x^*)$ (see Lemma 1.3). 
Together with the above Theorem D and Theorem 2.1 in [CRX], we obtain the 
following $C^{1,\alpha}$-compactness result.

\proclaim{Theorem E} Given $n, \rho, d, \Lambda, v>0$, there exist 
$\epsilon=\epsilon(n,\rho,d,\Lambda,v)$ such that the collection of compact 
$n$-manifolds satisfying
$$\Lambda\ge \op{Ric}_M\ge (n-1)H,\, d\ge \op{diam}(M),\,\op{vol}(M)\ge v,\, 
\frac{\op{vol}(B_\rho(x^*))}{\svolball{H}{\rho}}\ge 1-\epsilon, \,
\forall \, x\in M$$
is compact in the $C^{1,\alpha}$-topology, where the condition, `` 
$\op{vol}(M)\ge v$"  can be removed when $H=-1$.
\endproclaim

A few remarks are in order:

\remark{Remark \rm 0.5} In the proof of Theorem A, we actually proved that
$g$ is almost Einstein on any $B_{\frac \rho8}(x)$, $x\in M$
(see (2.1.2)); compare to Problem 2.4. Roughly, one may interpret this as
under bounded Ricci curvature, a ball with almost maximal `rewinding volume'
is an almost `Einstein ball'.
\endremark

\remark{Remark \rm 0.6} The existence of a nearby metric of almost constant
curvature in Theorem B is crucial to our proof of Theorem C (and Theorem D).
Indeed, we do not know, even assuming a higher regularity on
the original metric, how to prove Theorem C without using a nearby
metric of almost constant sectional curvature.
\endremark

\remark{Remark \rm 0.7} In Theorem C, no restriction on diameter for
$H=-1$. For $H=0$, the condition on bounded diameter cannot be removed.
Here is a counterexample: for each $i$, let $S^3_i$ denote a round $3$-sphere
of radius $i$, and let $g_i$ be a collapsed Berger's metric such that
$\op{vol}(B_1(p,g_i))<i^{-1}$ and $i^{-5}<\op{sec}(g_i)< 4i^{-2}$ (p.81, [Pet]).
It is easy to see that $\frac{\op{vol}(B_1(p^*,g^*_i))}{\svolball{0}{1}}\to 1$,
as $i\to \infty$.
\endremark

\remark{Remark \rm 0.8} Theorem D verifies Conjecture 0.2 (thus Conjecture 0.3)
for the case that Ricci curvature is bounded above. Note that for $H=0,1$, $M$ in
Theorem D can be collapsed (comparing to Theorem 0.4 where $M$ is not collapsed
when $H=1$).
\endremark

\remark{Remark \rm 0.9} Note that Theorem E and the $C^{1,\alpha}$-compactness
theorem in [An] may have only a `small' overlap. This is because
the local volume condition in Theorem E and the injectivity radius condition
in [An] are somewhat `parallel': a lower bound on injectivity radius
may not imply the volume condition in Theorem E,
and vice versa the volume conditions may not imply a lower bound
on injectivity radius. (note that for $H=-1$, in Theorem E  a priori $M$
could be collapsed.)
\endremark

\vskip4mm

The rest of the paper is organized as follows:

In Section 1, we will supply notions and basic properties that
will be used through out the rest of the paper.

In Section 2, we will prove Theorems A-E. At the end, we will
ask a few questions relating to the approach in this paper.

\vskip4mm

{\bf Acknowledgement}: The authors would like to thank Jian Song,
Zhenlei Zhang, and Bin Zhou for helpful discussion on Ricci flows.

\vskip8mm

\head 1. Preliminaries
\endhead

\vskip4mm

The purpose of this section is to supply notions and basic properties that
will be used through out the rest of the paper; we refer readers to [An], [CC1]
and [DWY] for details.

\vskip4mm

\subhead a. Almost maximal volume ball is an almost space form ball
\endsubhead

\vskip4mm

Let $N$ be a Riemannian $(n-1)$-manifold, let $k:(a,b)\to \Bbb R$ be a smooth positive
function and let $(a,b)\times_k N$ be the $k$-warped product whose Riemannian tensor is
$$g = dr^2 + k^2(r)g_N.$$
The Riemannian distance $|(r_1,x_1)(r_2,x_2)|$ ($x_1\neq x_2$) equals to the
infimum of the length
$$\int_0^l\sqrt{ (c_1'(t))^2+k^2(c_1(t))}dt$$
 for any smooth curve $c(t)=(c_1(t),c_2(t))$ such that $c(0)=(r_1,x_1)$, $c(l)=(r_2,x_2)$
 and $|c_2'|\equiv 1$, and $|(r_1,x)(r_2,x)|=|r_2-r_1|$.
Thus given $a,b,k$, there is a function (e.g., the law of cosine on space forms)
$$\rho_{a,b,k}(r_1,r_2,|x_1x_2|)=|(r_1,x_2)(r_2,x_2)|.$$
Using the same formula for $|(r_1,x_2)(r_2,x_2)|$,
one can extend the $k$-warped product $(a,b)\times_k Y$ to any metric space $Y$ (not necessarily
a length space); see \cite{CC1}.

The following theorem in [CC1] asserts that an almost volume annulus (see (1.1.1) below)
is an almost metric annulus (see (1.1.2)).

\proclaim{Theorem 1.1 (\cite{CC1})}
Let $M$ be a Riemannian manifold, let $r$ be a distance
function to a compact subset in $M$, let $A_{a,b}=r^{-1}((a,b))$, let
 $$\Cal V(u)=\inf\left\{\left.\frac{\volume (B_u(q))}{\volume (A_{a,b})}
 \right| \text{ for all $q\in A_{a,b}$ with $B_u(q)\subset A_{a,b}$ }\right\},$$
and let $0<\alpha'<\alpha, \alpha-\alpha'>\xi>0$.  If
$$\op{Ric}_M \ge -(n-1)\frac{k''(a)}{k(a)}\qquad (\text{on $r^{-1}(a)$}),$$
$$\Delta r\le (n-1)\frac{k'(a)}{k(a)} \qquad (\text{on $r^{-1}(a)$}),$$
$$\frac{\volume(A_{a,b})}{\volume(r^{-1}(a))}\ge (1-\epsilon)
\frac{\int_a^bk^{n-1}(r)dr}{k^{n-1}(a)}.\tag 1.1.1$$
Then there exists a length metric space $Y$, with at most $\#(a,b,k,\Cal{V})$
components $Y_i$, satisfying
$$\op{diam}(Y_i)\le D(a,b,k,\Cal V),$$
such that
the Gromov-Hausdorff distance,
$$d_{GH}(A_{a+\alpha,b-\alpha},(a+\alpha,b-\alpha)\times_kY)\le
\Psi(\epsilon|n,k,a,b,\alpha',\xi,\Cal V) \tag 1.1.2$$
with respect to the two metrics $d^{\alpha',\alpha}$ and $\b d^{\alpha',\alpha}$, where
$d^{\alpha',\alpha}$ (resp. $\b d^{\alpha',\alpha}$) denotes the restriction of the
intrinsic metric of $A_{a+\alpha',b-\alpha'}$ on $A_{a+\alpha,b-\alpha}$ (resp.
$(a+\alpha',b-\alpha')\times_kY)$ on $(a+\alpha,b-\alpha)\times_kY$).
\endproclaim

Let $$\op{sn}_H(r)=\cases \frac{\sin \sqrt Hr}{\sqrt H} & H>0\\ r & H=0\\
\frac{\sinh \sqrt{-H}r}{\sqrt{-H}} & H<0\endcases.$$
Applying Theorem 1.1 to $k=\op{sn}_H(r)$ with $r(x)=d(p,x): M\to \Bbb R$,
we conclude the following result that is used in the proof of Theorem A-E.

\proclaim{Theorem 1.2} For $n, \rho, \epsilon >0$, if a complete
$n$-manifold $M$ contains a point $p$ satisfies
$$\op{Ric}_M\ge (n-1)H,\quad \frac{\volume(B_\rho(p))}{\svolball{H}
{\rho}}\ge 1-\epsilon,$$
then $d_{GH}(B_{\frac \rho2}(p),\sball{H}{\frac \rho2})
<\Psi(\epsilon|n,\rho,H)$.
\endproclaim

Note that $\frac{\volume(B_\rho(p))}{\svolball{H} {\rho}}\ge 1-\epsilon$
implies (1.1.1), as $a\to 0$. Since the almost maximal
volume condition holds at all points near $p$ (which contains regular
points), by simple blow up argument one concludes that $Y$ is isometric to
$S^{n-1}_1$.
\vskip4mm

\subhead b.  Almost maximal volume and $C^{1,\alpha}$-harmonic radius estimate
\endsubhead

\vskip4mm

In this and the next subsections, we will always assume bounded Ricci curvature:
$\Lambda\ge \op{Ric}\ge (n-1)H$, $H=\pm 1$ or $0$.

Let $M$ be a complete $n$-manifold. For $p\in M$, $k\ge 0$, $0<\alpha<1$
and $Q\ge 1$,  the $C^{k,\alpha}$-harmonic radius at $p$
with respect to $Q$ is the largest radius $r_h(p)$ of the ball at $p$ such that there
are harmonic coordinates on $B_{r_h}(p)$ and $r^{k+\alpha}|g_{ij}|_{C^{k,\alpha},
B_{r_h}(p)}\le Q$. The harmonic radius of a subset is the infimum of the harmonic radii of
points in the subset.

\proclaim{Lemma 1.3}  For $n, \rho, \Lambda>0$, $Q>1$ and $0<\alpha<1$, there are
constants,  $\epsilon(n,\rho, \Lambda), r_h(n,\rho,\Lambda,Q,\alpha)>0$,  such that
if a compact Riemannian $n$-manifold $(M,g)$ satisfies
$$\Lambda\ge \op{Ric}(g)\ge (n-1)H,\quad \frac{\op{vol}(B_\rho(x^*))}{\svolball{H}
{\rho}}\ge 1-\epsilon(n,\rho,\Lambda), \quad \forall \, x\, \in M,$$
the $C^{1,\alpha}$-harmonic radius on $B_{\frac \rho2}(x^*)$ with respect to
$Q$ is at least $r_h(n,\rho,\Lambda,Q,\alpha)$.
\endproclaim

\demo{Proof}  We argue by contradiction, and the proof follows the same argument
as in the proof Main Lemma 2.2 in [An]; where the almost maximal volume condition is replaced with
a lower bound on injectivity radius which is to guarantee that any blow up limit
is $\Bbb R^n$. We claim that a contradicting sequence, $\frac{\op{vol}(B_\rho(x_i^*))}
{\svolball{H}{\rho}}\ge 1-\epsilon_i\to 1$, also satisfies that any blow up limit is $\Bbb R^n$.
Hence, the same proof in [An] goes through here to derive a contradiction.

To see the claim, for any $y_i^*\in B_{\frac \rho2}(x_i^*)$, $R>0$ and $r_i\to \infty$, by
Bishop-Gromov relative volume comparison and the volume convergence in [Co2] we derive
$$\frac{\op{vol}(B_R(y_i^*,r^2_ig_i^*))}{\svolball{r_i^{-2}H}R}=\frac{\op{vol}(B_{r_i^{-1}R}(y_i^*))}
{\svolball{H}{r_i^{-1}R}}\ge \frac{\op{vol}(B_{\frac \rho8}(y_i^*))}
{\svolball{H}{\frac\rho8}}\ge 1-\Psi(\epsilon_i|n,\rho).$$
Since $r_i^{-2}H\to 0$, by Theorem 1.2 we conclude that
$$d_{GH}(B_R(y_i^*,r^2_ig_i^*),\sball{0}{R})\to 0.$$
Since $R$ is arbitrarily chosen, the desired claim follows.
\qed\enddemo

As an application of Lemma 1.3, we will prove a non-collapsed local version of
Theorem A.

\proclaim{Lemma 1.4}  Given $n, \rho, \Lambda >0$ and $H=\pm 1$ and $0$, there is
$\epsilon(n,\rho,\Lambda)>0$ such that for $0<\epsilon<\epsilon(n,\rho,\Lambda)$, if a
complete $n$-manifold $(M,g,x)$ satisfies
$$\Lambda \ge \op{Ric}(g)\ge (n-1)H,\quad \frac{\volume(B_{\rho}(x))}
{\svolball{H}{\rho}}\geq 1-\epsilon,$$
then for all $p\ge 1$,
$$-\kern-1em\int_{B_{\frac \rho2}(x)}|\op{Ric}(g)-(n-1)H g|^p\leq \Psi(\epsilon |n,\rho,\Lambda,p).$$
\endproclaim

\demo{Proof} Arguing by contradiction, assume a contradicting sequence, $(M_i,g_i,x_i)$,
satisfying
$$\Lambda \ge \op{Ric}_{M_i}\ge (n-1)H,\quad \frac{\volume(B_{\rho}(x_i))}
{\svolball{H}{\rho}}\ge 1-\epsilon_i\to 1,$$
but $-\kern-1em\int_{B_{\frac \rho2}(x_i)}|\op{Ric}(g_i)-(n-1)H g_i|^{p_0}\ge \delta_0>0$, for some $p_0\ge 1$.

By Theorem 1.2, we may assume that $B_{\frac \rho2}(x_i)@>GH>>\sball{H}{\frac \rho2}$. By [CC2],
we may assume that for $i$ large, $B_{\frac \rho2}(x_i)$ are diffeomorphic to
$\sball{H}{\frac \rho2}$. From the expression of Ricci curvature
in a harmonic coordinate, a bound on Ricci curvature implies
that $g_i\to \b g_H$ in $L^{2,p}$-norm for all $p\ge 1$. Consequently,
$h_i=\op{Ric}(g_i)-(n-1)Hg_i\to h=\op{Ric}(g)-(n-1)Hg\equiv 0$ on $\sball{H}{\frac \rho2}$
in $L^p$-norm, a contradiction.
\qed\enddemo

\vskip4mm

\subhead c. Almost maximal volume and Ricci flows
\endsubhead

\vskip4mm

The main reference for this subsection is [DWY].

Let $(M,g)$ be a compact Riemannian manifold. The Ricci flow was introduced by
Hamilton as the solution of the following parabolic PDE,
$$\frac{\partial}{\partial t}g(t)=-2\op{Ric}(g(t)),\qquad g(0)=g.$$
The solution always exists for a short time $t>0$, and if the maximal flow time
$T_{\op{max}}<\infty$,
then $\max |\op{Rm}(g(t))|\to +\infty$ as $t\to T_{\op{max}}$ ([Ha]).

A basic property of Ricci flow is that it improves the regularity
of the initial metric ([Sh1,2]). However, the regularity depends on the
flow time. For our purpose, a uniform definite
flow time is important. We have

\proclaim{Theorem 1.5} For $n, \rho, \Lambda>0$ and $H=\pm 1$ or $0$, there are positive
constants, $\epsilon(n,\rho,\Lambda),  T(n,\rho,\Lambda)$, such that if a compact Riemannian
$n$-manifold $(M,g)$ satisfies
$$\Lambda\ge \op{Ric}(g)\ge (n-1)H,\quad \frac{\op{vol}(B_\rho(x^*))}{\svolball{H}{\rho}}
\ge 1-\epsilon(n,\rho,\Lambda), \quad \forall \, x\, \in M,$$
then the Ricci flow,
$$\frac{\partial }{\partial t}g(t)=-2\op{Ric}(g(t)),\quad g(0)=g,$$
exists for $t\in [0,T(n,\rho,\Lambda)]$ and
$$|g(t)-g|<4t,\quad |\op{Rm}(g(t))|_{C^k}\le C,
\quad \Lambda+ct^{\frac 12}\ge \op{Ric}(g(t))\ge (n-1)H-ct^{\frac 12},$$
where $C=C(n,\rho,\Lambda,k,t)$ and $c=c(n,\rho,\Lambda)$.
\endproclaim

Note that Theorem 1.5 is similar to Theorem 1.1 in [DWY], where the volume condition
on local covering is replaced by a positive lower bound on conjugate radius. Note that
the condition on conjugate radius is solely used to show a $L^{2,p}$-harmonic
radius lower bound on a local covering space for all $p\ge 1$ (see Remark 1 in [DWY]),
which is required to apply the Moser's weak Maximum principle 
(Theorem 2.1 in [DWY]).
Because a lower bound on the $L^{2,p}$-harmonic radius follows from Lemma 1.3 and
bounded Ricci curvature condition, the same proof in [DWY] will give a proof
of Theorem 1.5 with the obvious modification (cf. [Sh1,2]).

Let $(M,g)$ be as in Theorem 1.5. Inspired by [DWY] we will show
that if $g$ is almost $H$-Einstein in $L^p$-sense, then the renormalized
Ricci flow solution $g(t)$ in (1.6.1) below is again almost $H$-Einstein
in $L^p$-sense (Lemma 1.7).

Consider the renormalized Ricci flow in the sense of [TW]:
$$\frac{\partial}{\partial t} g=-\op{Ric}(g)+(n-1)Hg,\tag 1.6.1$$
and let
$$\bar g(s)=\cases
\sqrt{1-4(n-1)Hs}\cdot g\left
(\frac{\ln(1-4(n-1)Hs)}{-2(n-1)H}\right), &
H=\pm 1\\
g(2s), & H=0.
\endcases$$
Then $\bar g(s)$ satisfies $\bar g(0)=g(0)$ and
$$\frac{\partial }{\partial s}\bar g=-2\op{Ric}(\bar g(s)). \tag 1.6.2$$
Let $g^*(t)$ (resp. $\bar g^*(s)$) be the lifting of $g(t)$
(resp. $\bar g(s)$) on $B_\rho(x^*)$. Then
$$|\bar R^*_{ijkl}(s)|=\cases
\frac 1{\sqrt{1-4(n-1)Hs}}\cdot R^*_{ijkl}\left
(\frac{\ln(1-4(n-1)Hs)}{-2(n-1)H}\right), &
H=\pm 1\\
R^*_{ijkl}(2s), & H=0.
\endcases$$
Let
$$h^*_{ij}=R^*_{ij}-(n-1)H g^*_{ij}.\tag 1.6.3$$
Then
$$\frac{\partial }{\partial t}h^*_{ij}=\frac{1}{2}\Delta
h^*_{ij}+R^*_{pijq}h^*_{pq}-h^*_{ip}h^*_{pj}.$$
To get ride of the $\frac 12$-factor, we make a change of variable $t=2t'$
(for simple notation, switch by $t'=t$). Then the above implies
$$\frac{\partial }{\partial t}|h^*|\le \Delta
|h^*|+2|Rm^*||h^*|.\tag 1.6.4$$
By applying the Moser's weak maximum principle, we
conclude the following:

\proclaim{Lemma 1.7} Let the assumptions be in Theorem 1.5,
and let $h_{ij}^*(t)$ be defined in the above, and
$$\bar T(n,\rho,\Lambda)=\cases
\frac{\ln(1-4(n-1)HT(n,\rho,\Lambda))}{-4(n-1)H}, &
H=\pm 1\\
T(n,\rho,\Lambda), & H=0.
\endcases$$
Then for $t\in (0,\bar T(n,\rho,\Lambda)]$,
$$\max_{x\in M}|h^*(t)|_{p,B_{\frac\rho4}(x^*,g)}\le \max_{x\in
M}|h^*(0)|_{p,
B_{\frac\rho2}(x^*,g)}\cdot
\frac{1}{1-c(n,\rho,\Lambda)t}.$$
\endproclaim

\demo{Proof} Given Lemma 1.3, by (1.6.4) the rest proof is an imitation of
the proof of Lemma 3.3 in [DWY].
\qed\enddemo

\vskip4mm

\head 2. Proof of Theorems A-E
\endhead

\vskip4mm

\demo{Proof of Theorem A}

Because Ricci curvature is bounded in absolute value, it suffices to prove Theorem A for $p=1$.
For any $x\in M$, by Lemma 1.4 we have
$$-\kern-1em\int_{B_{\frac \rho2}(x^*)}|\op{Ric}(g^*)-(n-1)H g^*|\leq
\Psi(\epsilon | n,\rho,\Lambda).\tag 2.1.1$$
We claim that for $\Psi_1(\epsilon|n,\rho,\Lambda)=\Psi(\epsilon|n,
\rho,\Lambda)\cdot \frac{\svolball{H}{\frac \rho2}}{\svolball{H}{\frac \rho8}}$,
the following holds:
$$-\kern-1em\int_{B_{\frac\rho8}(x)}|\op{Ric}(g)-(n-1)H g|\leq
\Psi_1(\epsilon |n,\rho,\Lambda).\tag 2.1.2$$
Let $A=\{x_i\}$ denote an $\frac \rho8$-net on $M$. Then
$B_{\frac\rho{16}}(x_i)\cap B_{\frac\rho{16}}(x_j)=\emptyset$ ($i\ne j$)
and $M\subseteq \bigcup_{x_i\in A}B_{\frac\rho8}(x_i)$.
Assuming (2.1.2), we derive
$$\align
-\kern-1em\int_M |\op{Ric}(g)-(n-1)H g|&\le \frac{1}{\volume(M)}\sum_{x_i\in A}
\int_{B_{\frac\rho8}(x_i)}|\op{Ric}(g)-(n-1)H g|\\
&= \frac{1}{\volume(M)}\sum_{x_i}\volume(B_{\frac\rho8}(x_i))-
\kern-1em\int_{B_{\frac\rho8}(x_i)}|
\op{Ric}(g)-(n-1)H g|\\ &\leq \frac 1{\op{vol}(M)}\sum_{x_i\in A}
\volume(B_{\frac\rho8}(x_i))\Psi_1(\epsilon | n,\rho,\Lambda)\\
&\leq \frac 1{\op{vol}(M)}\sum_{x_i\in A} \volume(B_{\frac\rho{16}}(x_i))\cdot
\frac{\svolball{H}{\frac\rho8}}{\svolball{H}{\frac\rho{16}}}\Psi_1(\epsilon | n,\rho,\Lambda)\\
&\leq \Psi_2(\epsilon | n,\rho,\Lambda).
\endalign$$

We now verify (2.1.2). Let $D$ denote the Dirichlet fundamental domain at $x^*\in
\widetilde {B_\rho(x)}$, and let $\Gamma(\frac\rho4) =
\{\gamma\in\pi_1(B_{\rho}(x^*)),\,\;|x^*\gamma (x^*)|\le \frac\rho4\}$.
Then
$$B_{\frac\rho8}(x^*)\subset \bigcup_{\gamma\in \Gamma(\frac\rho4)}
\gamma(B_{\frac\rho8}(x^*)\cap D)\subset  B_{\frac \rho2}(x^*),$$
We claim that there is a $\gamma\in \Gamma(\frac\rho4)$ such that
$$-\kern-1em\int_{\gamma(B_{\frac \rho8}(x^*)\cap D)}|\op{Ric}(g^*)-(n-1)H g^*|
\leq \Psi_1(\epsilon | n,\rho,\Lambda),$$
i.e.,
$$-\kern-1em\int_{B_{\frac\rho8}(x)}|\op{Ric}(g)-(n-1)H g|\leq \Psi_1(\epsilon | n,\rho,\Lambda).$$
If the claim fails i.e., for all $\gamma\in \Gamma(\frac\rho4)$,
$$-\kern-1em\int_{\gamma(B_{\frac\rho8}(x^*)\cap D)}|\op{Ric}(g^*)-(n-1)H g^*|>
\Psi_1(\epsilon | n,\rho,\Lambda),$$
then
$$\align -\kern-1em\int_{B_{\frac \rho2}(x^*)}&|\op{Ric}(g^*)-(n-1)H g^*|\\ &\ge
\sum_{\gamma\in \Gamma(\frac \rho4)}\frac {\op{vol}(\gamma(B_{\frac\rho8}(x^*)\cap D))}
{\op{vol}(B_{\frac \rho2}(x^*))}-\kern-1em
\int_{\gamma(B_{\frac\rho8}(x^*)\cap D)}|\op{Ric}(g^*)-(n-1)H g^*|\\
&> \frac {\Psi_1(\epsilon|n,\rho,\Lambda)}{\op{vol}(B_{\frac \rho2}(x^*))}\sum_{\gamma\in \Gamma(\frac\rho4)}
\volume(\gamma(B_{\frac\rho8}(x^*)\cap D))\\
&\geq \frac{\Psi_1(\epsilon | n,\rho,\Lambda)\volume(B_{\frac\rho8}(x^*))}{\volume(B_{\frac \rho2}(x^*))}\\&
\ge \Psi(\epsilon|n,\rho,\Lambda),
\endalign$$
a contradiction to (2.1.1).
\qed \enddemo

\demo{Proof of Theorem B}

Arguing by contradiction, assume a contradicting sequence, $(M_i,g_i)@>GH>>X$,
such that
$$\Lambda\ge \op{Ric}(g_i)\ge (n-1)H,\quad \frac
{\op{vol}(B_\rho(x_i^*))}{\svolball{H}{\rho}}
\ge 1-\epsilon_i\to 0,\, \forall\, x_i\in M_i,$$
and $M_i$ admits no nearby metric to $g_i$ with almost constant sectional curvature $H$.

Fixing a small $\delta\in (0,T(n,\rho,\Lambda)]$ (Theorem 1.5),
let $g_i(\delta)$ denote the renormalized Ricci flow in (1.6.1).
By Theorem 1.5, for any $x_i\in M_i$, passing to a subsequence we may assume
that the lifting metric $g_i^*(\delta)$ on $B_\rho(x_i^*)$ satisfies
$$B_{\frac \rho2}(x^*_i,g_i^*(\delta))@>C^k>>B_{\frac \rho2}(x_\delta^*,g^*_\infty(\delta)),\quad
h_i(g_i^*(\delta))@>C^k>> h(g^*_\infty(\delta)),$$
where $h_i$ is defined in (1.6.3), and the $C^k$-convergence can be seen from
the Cheeger-Gromov convergence theorem. Consequently, $g^*_\infty(\delta)$
is a smooth metric and $h(g_\infty^*(\delta))$ is a
smooth tensor on $B_{\frac \rho2}(x_\delta^*,g^*_\infty(\delta))$. By Lemma 1.4 and
Lemma 1.7, for any $x_i\in M_i$,
$$|h(g^*_i(\delta))|_{p,B_{\frac \rho4}(x_i^*)}\to 0.$$
Consequently, $h(g_\infty^*(\delta))|_{B_{\frac \rho4}(x^*_\delta,g_\infty^*(\delta))}\equiv 0$ i.e.,
$g^*_\infty(\delta)|_{B_{\frac \rho4}(x^*_\delta,g_\infty^*(\delta))}$ is $H$-Einstein.

Clearly, $B_{\frac \rho4}(x_\delta^*,g^*_\infty(\delta))@>GH>>\sball{H}{\frac \rho4}$ as $\delta\to 0$.
Since $g_\infty^*(\delta)$ is $H$-Einstein for all $\delta$,
$B_{\frac \rho4}(x^*_\delta,g^*_\infty(\delta))@>C^k>>\sball{H}{\frac \rho4}$,
for any $k\ge 1$ ([CC2]). Consequently,
for $\delta_0$ sufficiently small, $g^*_\infty(\delta_0)$ has almost constant sectional curvature $H$.
Since $B_{\frac \rho4}(x_i^*,g_i^*(\delta_0))@>C^k>>B_{\frac \rho4}(x_{\delta_0}^*,g^*_\infty(\delta_0))$,
for $i$ large, $g_i^*(\delta_0)$ has almost constant curvature $H$.
Since $x_i$ is arbitrarily chosen, we conclude that $g_i(\delta_0)$ has almost
constant sectional curvature $H$, a contradiction.
\qed\enddemo

\demo{Proof of Theorem C}

Fixing a small $\delta>0$, by Theorem B we may assume a nearby metric $g(\delta)$
such that
$$|g-g(\delta)|<\delta,\quad H-\delta\le \op{sec}_{g(\delta)}\le H+\delta.$$

Case 1. Assume $H=-1$. For any $\tilde p\in \tilde M$, the exponential map,
$\exp^{\tilde g(\delta)}_{\tilde p}: T_{\tilde p}\tilde M\to \tilde M$, is
a diffeomorphism such that its differential has a bounded norm on $B_1(0)$
depending on $n$. Consequently, $\op{vol}(B_1(\tilde p,\tilde g(\delta))$ has a positive
lower bound depending only on $n, \rho$ and $d$. Since $|\tilde g-
\tilde g(\delta)|<\delta$, we conclude the desired result.

Case 2. Assume $H=0$. By Splitting theorem of Cheeger-Gromoll,
$\tilde M=\Bbb R^k\times N$, where $N$ is a simply connected
$(n-k)$-manifold of non-negative Ricci curvature. We claim that
$N$ is a point. Note that $\op{diam}(N)\le c(n,d)$ (see the proof
of Theorem C in [CRX] where we normalize $d=1$). We may assume $\delta^{-\frac 12}>4\op{diam}(N)$. Note that
since $\op{sec}_{g(\delta)}<\delta$, $\exp^{\tilde g(\delta)}_{\tilde p}: B_{\frac 1{\sqrt\delta}}(0)\to
B_{\frac 1{\sqrt\delta}}(\tilde p,\tilde g(\delta))$ is a local diffeomorphism.
Note that $B_{\frac 1{2\sqrt\delta}}(\tilde p)$ can be deformed to $0\times N$ ($\tilde p=(0,x)$)
and thus $B_{\frac 1{2\sqrt \delta}}(\tilde p)$ is simply connected. Consequently, the
lifting of $B_{\frac 1{2\sqrt\delta}}(\tilde p)$ via $\exp_{\tilde p}^{\tilde g(\delta)}$ is contained
in the segment domain (i.e. each $\tilde x\in B_{\frac 1{2\sqrt \delta}}(\tilde p)$ is connecting
to $\tilde p$ by a unique minimal geodesic; if $c_1$ and $c_2$ are two distinct minimal geodesics,
then $c_1*c^{-1}_2$ is a loop at $\tilde p$, and so is the lifting of $c_1*c^{-1}_2$
a loop at $0$. Note that with respect to the pullback metric on $T_{\tilde p}\tilde M$,
we obtain two geodesics from $0$ to some $v$; a contradiction).
Therefore, $B_{\frac 1{2\sqrt\delta}}(\tilde p)$ is contractible
in $\tilde M$, a contradiction.

Case 3. Assume $H=1$. The classical $1/4$-pinched injectivity radius estimate
implies that the pullback metric $\tilde g(\delta)$ on $\tilde M$ has injectivity radius $>\frac \pi2$, and thus
$\op{vol}(B_1(\tilde p))$ has a positive lower bound depending on $n$.
By now the desired result follows.
\qed\enddemo

In the rest of the paper, we will freely use properties of equivariant Gromov-Hausdorff
convergence; see b. of Section 1 in [CRX] for details.

\proclaim{Lemma 2.2}  Let $M_i$ be a sequence of compact $n$-manifolds
satisfying
$$\op{Ric}_{M_i}\ge (n-1),\quad |\op{Rm}|_{C^1,M_i}\le C,
\quad \op{vol}(\tilde M_i)\ge v>0,$$
and the commutative diagram,
$$\CD (\tilde M_i,\Gamma_i)@>GH>>(\tilde M_\infty,G)
\\@VV \pi_i V   @VV \pi  V \\
M_i@>GH>>X,
\endCD$$
Then for $i$ large,

\noindent (2.2.1) There is injective homomorphism and
$\epsilon_i$-GHA ($\epsilon_i\to 0$), $\phi_i: \Gamma_i\to G$, such that
$\phi_i(\Gamma_i)$ acts freely on $\tilde M_\infty$.

\noindent (2.2.2) There is a $\Gamma_i$-conjugate
diffeomorphism, $\tilde f_i: (\tilde M_i,\Gamma_i)
\to (\tilde M_\infty,\phi_i(\Gamma_i))$, which is also an $\epsilon_i$-GHA.
\endproclaim

\demo{Proof}  Lemma 2.2 is essentially Theorem 3.5 in [CRX] where condition,
``$|\op{Rm}|_{C^1,M_i}\le C$'' , is replaced with ``$\frac{\op{vol}(B_\rho(\tilde x_i))}
{\svolball{1}{\rho}}\ge 1-\epsilon_i\to 1$'', and the proof of Theorem 3.5
proves Lemma 2.2 with the following modifications: the regularity condition
in Lemma 2.2 implies the following:

(i) There is a uniform lower bound on the injectivity radius of $\tilde M_i$, and thus
$\tilde M_\infty$ is a Riemannian manifold,  and for any $r_i\to \infty$, passing to a
 subsequence $(\tilde M_i,\tilde p_i,r_i^2\tilde g_i)$
converges to $\Bbb R^n$; which guarantees (2.2.1).

(ii) $\tilde f$ in (2.2.2) is a diffeomorphism, instead of a homotopy equivalence in Theorem 3.5;
see the discussion following Theorem 3.5.
\qed\enddemo

\demo{Proof of Theorem D}

By Theorem C, we may assume $\op{vol}(B_1(\tilde p))\ge v$. For $H\ne 1$,
by Theorem 0.4 we obtain the desired conclusion (indeed, the case $H=0$ has been
already proved in the proof of Theorem C).

For $H=1$, because $\op{vol}(M)$ can be very small, Theorem 0.4 cannot be applied
here. Arguing by contradiction, assume a
contradicting sequence, $(M_i,g_i)$, such that $g_i$ satisfies
the conditions of Theorem D for $\epsilon_i\to 0$ but none of $M_i$ is
diffeomorphic to a spherical space form.

For each $i$, let $g_i(\delta)$ be as in Theorem B, such that for all $1\le k<\infty$,
$$|\op{Rm}(g_i(\delta))|_{C^k}\le C(n,\rho,\Lambda,\delta,k),\quad
1-\delta\le\op{sec}_{g_i(\delta)}\le 1+\delta.$$
Passing to a subsequence we may assume
the following commutative diagram:
$$\CD (\tilde M_i,\tilde g_i(\delta),\Gamma_i)@>GH>>(\tilde M_\infty(\delta),\tilde g_\infty(\delta),G(\delta))
\\@VV \pi_i V   @VV \pi  V \\
(M_i,g_i(\delta))@>GH>>(X,d_\infty(\delta)),
\endCD$$
where $\Gamma_i$ denotes the deck transformations. Since $\tilde M_i$ is not collapsed
(Theorem C), by Lemma 2.2 there is a $\Gamma_i$-conjugate diffeomorphism,
$\tilde f_i(\delta): (\tilde M_i,
\tilde g_i(\delta),\Gamma_i)\to (\tilde M_\infty(\delta),\tilde g_\infty(\delta),
\phi_i(\delta)(\Gamma_i))$. From the proof of Theorem B, we see that
$(\tilde M_\infty(\delta),\tilde g_\infty(\delta))$ is $1$-Einstein. It is clear that
$(\tilde M_\infty(\delta),\tilde g_\infty(\delta),G(\delta))@>GH>>(S^n_1,\b g_1,G)$, as $\delta\to 0$.
Consequently, for all $k<\infty$, $(\tilde M_\infty(\delta),\tilde g_\infty(\delta),G(\delta))
@>C^k>>(S^n_1,\b g_1,G)$ ([CC2]).

For each $\delta$, we may choose $i$ large such that $d_{GH}(\phi_i(\delta)(\Gamma_i),G(\delta))<\delta_i\to 0$
i.e., $(M_\infty(\delta),\tilde g_\infty(\delta),\phi_i(\delta)(\Gamma_i))@>GH>>(S^n_1,\b g_1,G)$.
We then apply Lemma 2.2 again to conclude that for a fixed small $\delta$,
there is $\phi_i(\delta)(\Gamma_i)$-conjugate diffeomorphism, $\tilde f_\infty(\delta):
(\tilde M_\infty(\delta),\phi_i(\delta)(\Gamma_i))\to (S^n_1,\psi_i(\delta)\circ \phi_i(\delta)(\Gamma_i))$. Then
$\tilde f_\infty(\delta)\circ \tilde f_i(\delta): (\tilde M_i,\Gamma_i)\to (S^n_1,\psi_i(\delta)\circ \phi_i(\delta)(\Gamma_i))$ is
$\Gamma_i$-conjugate diffeomorphism, and thus $M_i$ is diffeomorphic to a spherical space form,
$S^n_1/(\psi_i(\delta)\circ \phi_i(\delta)(\Gamma_i))$,
a contradiction.
\qed\enddemo

\remark{Remark \rm 2.3} Given Theorem B, the conclusion of Theorem D for $H=0$ and $H=1$
can also be seen from the work [Gr] and [BS] respectively.
\endremark

\demo{Proof of Theorem E}

It suffices to show that for any $Q\ge 1$ and $0<\alpha<1$, there is a constant $r_h=r_h(n,\rho,d,\Lambda,v,\alpha,Q)>0$ such that $M$ has $C^{1,\alpha}$-harmonic radius with respect to $Q$ bounded below by $r_h$; because $\Lambda\ge \op{Ric}_M\ge (n-1)H$.

Arguing by contradiction, assume  for some $Q_0\ge 1$ and $0<\alpha_0<1$,  there is a contradicting sequence,
$M_i$, satisfying
$$\Lambda\ge \op{Ric}_{M_i}\ge (n-1)H,\quad d\ge \op{diam}(M_i),\quad \frac{\op{vol}(B_\rho(x_i^*))}
{\svolball{H}{\rho}}\ge 1-\epsilon_i\to 1,\quad \forall\, x_i\in M_i,$$
and $p_i\in M_i$ such that the $C^{1,\alpha_0}$-harmonic radius $r_h(p_i)\to 0$. Passing to a subsequence, we may assume the following
commutative diagram:
$$\CD (\widetilde{B_\rho(p_i)},p^*_i,K_i)@>GH>>(X^*,p^*,K)
\\@VV \pi_i V   @VV \pi  V \\
(B_\rho(p_i),p_i)@>GH>>(B_\rho(p),p),
\endCD$$
where $K_i$ denotes the fundamental group of $B_\rho(p_i)$. Since
$\frac{\op{vol}(B_\rho(p^*_i))}{\svolball{H}{\rho}}\ge 1-\epsilon_i\to 1$,
by Theorem 1.2 we see that $B_\rho(p^*)$ is local isometric to a $H$-space form.
If $H\ne -1$, $K_i$ is discrete because $\op{vol}(M_i)\ge v$. We claim
that $K$ is discrete when $H=-1$. Hence, in any case we are able to
apply Theorem 2.1 in [CRX] to conclude that $K$ acts freely on $X^*$.
We may assume that any element in $K_i$
moves any $x_i^*$ in $B_{\frac \rho2}(p_i^*)$ at least $\delta$-distance,
where $\delta$ depends on $(X^*,K)$.
By Lemma 1.3, we may assume that $r_h(p^*_i)\ge r_h(n,\rho,\Lambda,\alpha_0,Q_0)>0$, and thus
$2r_h(p_i)\ge \min\{\delta, r_h(p^*_i)\}>0$, a contradiction.

To see that $K$ is discrete, note that by Theorem D we conclude that
$M_i$ is $\Psi(\epsilon|n,\rho,d,\Lambda)$ close to a hyperbolic manifold
$\Bbb H^n/\Gamma_i$. By Margulis-Heintze lemma ([He]), $\Bbb H^n/\Gamma_i$ is
not collapsed, and by the volume convergence in [Co2] we then conclude that
$M_i$ is not collapsed (so $B_\rho(x_i)$ is not collapsed), and thus $K$ is discrete.
\qed\enddemo

We will conclude this paper with the following questions related to the present
approach to Conjecture 0.3:

\example{Problem 2.4} Does Theorem A hold without an upper bound on Ricci curvature?
Indeed, it seems that even it is not known whether the scalar curvature is almost constant
in $L^p$-sense.
\endexample

\example{Problem 2.5} (Ricci flow time) For $n, \rho>0$, and $H=\pm 1$ or $0$,
are there constants, $\epsilon(n,\rho)>0,T(n,\rho)>0$,
such that for any $0<\epsilon<\epsilon(n,\rho)$, if a compact $n$-manifold $(M,g)$ satisfies
$$\op{Ric}(g)\ge (n-1)H,\quad \frac{\volume
(B_\rho(x^*))}{\svolball{H}{\rho}}\ge 1-\epsilon,\quad\forall\, x\in M,$$
then the Ricci flow from $g$ exists for $t\in [0,T(n,\rho)]$?
\endexample

\example{Problem 2.6} (Flows preserving almost Einstein) Let $(M,g)$ be a compact
$n$-manifold of $\op{Ric}_M\ge (n-1)H$ and
$$-\kern-1em\int_{B_{\frac \rho2}(p^*)}|\op{Ric}(g^*)-(n-1)Hg^*|<\epsilon.$$
Let $g(t)$ be a renormalized Ricci flow of $g$ (see (1.6.1)). Is $-\kern-1em\int_{B_{\frac \rho2}
( p^*,g^*(t))} |\op{Ric}(g^*(t))-(n-1)Hg^*(t)|<\Psi(\epsilon|n,\rho,t)$.
\endexample

Note that if there are affirmative answers to Problem 2.4-2.6, then
the approach in this paper can be extended toward a proof of Conjecture 0.2.


\vskip20mm

\Refs
\nofrills{References}
\widestnumber\key{APS1}

\vskip3mm

\ref
\key An
\by M. Anderson
\pages 429-445
\paper Convergence and rigidity of manifolds under Ricci
curvature bounds
\jour Invent. Math.
\vol 102
\yr 1990
\endref



\ref
\key BS
\by S. Brendle; R. Schoren
\pages 287-307
\paper Manifolds with 1/4-pinched curvature are space forms
\jour J.A.M.S.
\vol 22
\yr 2009
\endref

\ref
\key CC1
\by J. Cheeger; T. Colding
\pages 189-237
\paper Lower Bounds on Ricci Curvature and the Almost Rigidity of Warped Products
\jour Ann. of Math.
\vol 144, No. 1
\yr Jul., 1996
\endref

\ref
\key CC2
\by J. Cheeger; T. Colding
\pages 406-480
\paper On the structure of space with Ricci curvature bounded below I
\jour J. Diff. Geom
\vol 46
\yr 1997
\endref

\ref
\key Co1
\by T.H. Colding
\pages 193-214
\paper Large manifolds with positive Ricci curvature
\jour Invent. Math.
\vol 124 (1-3)
\yr 1996
\endref

\ref
\key Co2
\by T. Colding
\pages 477-501
\paper Ricci curvature and volume convergence
\jour Ann. of Math
\vol 145(3)
\yr 1997
\endref

\ref
\key CRX
\by L. Chen, X. Rong; S. Xu
\pages
\paper Quantitive volume rigidity of space form with lower Ricci curvature bound I
\jour Submitted
\vol
\yr
\endref

\ref
\key DWY
\by X. Dai, G. Wei; R. Ye
\pages 49-61
\paper Smoothing Riemannian metrics with Ricci curvature bounds
\jour Manuscrlpta Math
\yr 1996
\vol 90
\endref

\ref
\key Ha
\by R. Hamilton
\pages 255-306
\paper Three-manifolds with positive Ricci curvature
\jour J. Diff. Geom
\yr 1982
\vol 17
\endref

\ref
\key Gr
\by M. Gromov
\pages 231-241
\paper Almost flat manifolds
\jour J. Diff. Geom.
\vol 13
\yr 1978
\endref

\ref
\key He
\by E. Heintze
\pages
\paper Manningfaltigkeiten negativer Kriimmung
\jour Ph.D. thesis, Universit\'it Bonn Habilitationsschrift
\vol
\yr 1976
\endref

\ref
\key LW
\by F. Ledrappier; X. Wang
\pages 461-477
\paper An integral formula for the volume entropy with application to rigidity
\jour J. Diff. Geom.
\vol 85
\yr 2010
\endref

\ref
\key Pe
\by G. Perelman
\pages 299-305
\paper Manifolds  of  Positive  Ricci Curvature
with  Almost  Maximal  Volume
\jour J. AMS.
\vol 7,  No.  2
\yr Apr., 1994
\endref

\ref
\key Pet
\by P. Petersen
\pages
\paper Riemannian Geometry (second edition)
\jour Springer-Verlag, New York
\vol
\yr 2006
\endref

\ref
\key Sh1
\by W. Shi
\pages 223-301
\paper Deforming the metric on complete Riemannian manifolds
\jour J. Diff. Geom.
\vol 30
\yr 1989
\endref

\ref
\key Sh2
\by W. Shi
\pages 303-394
\paper Ricci deformation of the metric on complete non-compact Riemannian manifolds
\jour J. Diff. Geom.
\vol 30
\yr 1989
\endref

\ref
\key TW
\by G. Tian; B. Wang
\pages 1169-1209
\paper On the structure of almost Einstein manifolds
\jour J.A.M.S
\vol 28
\yr 2015
\endref


\endRefs
\enddocument